\documentclass[leqno, 12pt]{article}
\usepackage{amsthm}
\usepackage{color}
\usepackage{amsfonts}
\usepackage{amsmath}
\usepackage[margin=2.5cm]{geometry}
\usepackage{graphicx}
\usepackage{wrapfig}
\usepackage{amscd}
\usepackage{amssymb}
\usepackage{amstext}
\usepackage{lscape}
\usepackage{multirow}
\usepackage{caption}
\usepackage{lipsum}

\usepackage{array}
\title{\bf  Inegalit\u a\c ti de tip Chebyshev-Gr\" uss pentru operatorii Bernstein-Euler-Jacobi }
\author{Heiner Gonska${}^{*}$, Maria-Daniela Rusu${}^{*}$, Elena-Dorina St\u{a}nil\u{a}${}^{*}$}

\date{}

\newtheorem{theorem}{Teorem\u a}[section]
\theoremstyle{plain}

\newtheorem{definition}[theorem]{Defini\c{t}ie}

\newtheorem{lemma}[theorem]{Lem\u a}

\newtheorem{remark}{Remarc\u{a}}[section]
\numberwithin{equation}{section}

\begin{document}

\newcolumntype{L}[1]{>{\raggedright\let\newline\\\arraybackslash\hspace{0pt}}m{#1}}
\newcolumntype{C}[1]{>{\centering\let\newline\\\arraybackslash\hspace{0pt}}m{#1}}
\newcolumntype{R}[1]{>{\raggedleft\let\newline\\\arraybackslash\hspace{0pt}}m{#1}}

\newcommand*{\N}{\ensuremath{\mathbb{N}}}               
\newcommand*{\R}{\ensuremath{\mathbb{R}}}               
\newcommand*{\Q}{\ensuremath{\mathbb{Q}}}     

\newcommand\blfootnote[1]{%
  \begingroup
  \renewcommand\thefootnote{}\footnote{#1}%
  \addtocounter{footnote}{-1}%
  \endgroup
}

\maketitle

\begin{abstract}
	The classical form of Gr\"{u}ss' inequality was first published by
G. Gr\"{u}ss  and gives an estimate of the difference between the
integral of the product and the product of the integrals of two
functions. In the subsequent years, many variants of this inequality
appeared in the literature. The aim of this paper is to consider
some Chebyshev-Gr\"{u}ss-type inequalities and apply them to the Bernstein-Euler-Jacobi (BEJ) operators of first and second
kind. The first and second moments of the operators will be of great interest. 
\end{abstract}



\blfootnote{${}^{*}$University of Duisburg-Essen, Faculty of Mathematics, Forsthausweg 2, 47057 Duisburg, Germany\\
E-mails: heiner.gonska@uni-due.de, maria.rusu@uni-due.de, elena.stanila@stud.uni-due.de}
\section{Introducere}

\^{I}n cele ce urmeaz\u{a} vom prezenta rezultate clasice din literatur\u{a}. 

Inegalit\u{a}\c{t}ile de tip Chebyshev-Gr\"{u}ss au fost intens studiate de-a lungul anilor, mai ales datorit\u{a} numeroaselor aplica\c{t}ii. Ele reprezint\u{a} leg\u{a}tura \^{\i}ntre func\c{t}ionala lui Chebyshev \c{s}i inegalitatea lui G. Gr\"{u}ss. 

Func\c{t}ionala descris\u{a} de
	\begin{equation*}
		T(f,g):=\frac{1}{b-a}\int^{b}_{a}{f(x)g(x)\mathrm{d}x}-\frac{1}{b-a}\int^{b}_{a}{f(x)\mathrm{d}x}\cdot \frac{1}{b-a}\int^{b}_{a}{g(x)\mathrm{d}x},
	\end{equation*}
unde $f,g:[a,b]\to \R$ sunt func\c{t}ii integrabile, este binecunoscut\u{a} \^{\i}n literatur\u a drept func\c{t}ionala clasic\u{a} a lui Chebyshev. Pentru detalii, articolul  \cite{chebyshev:1907} poate fi de ajutor. 

Un prim rezultat pe care \^{i}l readucem \^ in aten\c{t}ia cititorului este dat de urm\u{a}toarea teorem\u a.

\begin{theorem} (vezi \cite{mitrinovici:1993})
	Fie $f,g:[a,b]\to\R$ dou\u{a} func\c{t}ii m\u{a}rginite \c{s}i integrabile, ambele cresc\u{a}toare sau descresc\u{a}toare. Mai mult, fie $p:[a,b]\to \R^{+}_0$ o func\c{t}ie m\u{a}rginit\u{a} \c{s}i integrabil\u{a}. Atunci avem
		\begin{equation}\label{ineg1}
			\int^{b}_{a}{p(x)\mathrm{d}x}\int^{b}_{a}{p(x)\cdot f(x)\cdot g(x)\mathrm{d}x}\geq \int^{b}_{a}{p(x)\cdot f(x)\mathrm{d}x}\int^{b}_{a}{p(x)\cdot g(x)\mathrm{d}x}.
		\end{equation}
	Dac\u{a} una din func\c{t}iile $f$ sau $g$ este cresc\u{a}toare \c{s}i cealalt\u{a} este descresc\u{a}toare, atunci inegalitatea (\ref{ineg1}) se inverseaz\u{a}. 
\end{theorem}

\begin{remark}
	Rela\c{t}ia (\ref{ineg1}) a fost introdus\u{a} pentru prima oar\u{a} de c\u{a}tre P. L. Chebyshev \^{i}n 1882 ( vezi \cite{chebyshev:1882}). Din acest motiv, este cunoscut\u{a} sub numele de inegalitatea lui Chebyshev. 
\end{remark}

\^{I}n continuare amintim unul din rezultatele esen\c{t}iale pe care se bazeaz\u a cercetarea noastr\u{a}, \c{s}i anume inegalitatea lui Gr\"{u}ss pentru func\c{t}ionala lui Chebyshev.

\begin{theorem}(vezi \cite{gruess:1935})
	Fie $f, g$ func\c{t}ii integrabile, definite pe intervalul $[a,b]$ cu valori \^{\i}n $\R$, astfel \^{i}nc\^{a}t $m\leq f(x)\leq M$, $p\leq g(x)\leq P$, pentru orice $x\in [a,b]$, unde $m,M,p,P\in \R$. Atunci are loc inegalitatea
		\begin{equation}
			\left|T(f,g)\right|\leq \frac{1}{4} (M-m)(P-p).
		\end{equation}
\end{theorem}

Urm\u{a}toarele inegalit\u{a}\c{t}i de tip Chebyshev-Gr\"{u}ss vor fi folosite \^{\i}n continuare, pentru a introduce rezultatele noastre.

\begin{theorem}(vezi \cite{rusu:2011})
 Dac\u{a} $f,g\in C[a,b]$ \c{s}i $x\in [a,b]$ este fixat, atunci are loc inegalitatea
 	\begin{equation}\label{inegCG1}
 		\left|T(f,g;x)\right|\leq \frac{1}{4}\widetilde{\omega}\left(f;2\cdot \sqrt{H((e_1-x)^2;x)}\right)\cdot \widetilde{\omega}\left(g;2\cdot \sqrt{H((e_1-x)^2;x)}\right).
 	\end{equation}
\end{theorem}

\begin{remark}
	Rezultatul de mai sus implic\u{a} folosirea celui mai mic majorant concav $\widetilde{\omega}$ al primului modul de continuitate. O defini\c{t}ie \c{s}i detalii cu privire la acesta se g\u{a}sesc, de exemplu, \^{i}n \cite{goraru:2013}. 
\end{remark}

\begin{remark}
Scopul nostru este s\u{a} aplic\u{a}m inegalitatea de mai sus \^{\i}n cazul operatorilor BEJ de tipul $I$ \c{s}i $II$, pentru diferite cazuri, lu\^{a}nd \^{\i}n considerare diferitele momente de ordinul doi. 
\end{remark}

\^{I}n cazul operatorilor liniari \c{s}i pozitivi care reproduc func\c{t}iile constante, dar nu le reproduc pe cele liniare, avem urm\u{a}torul rezultat.

\begin{theorem}(vezi Corolarul 5.1. din \cite{goraru:2013})
	Dac\u{a} $H:C[a,b]\to C[a,b]$ este un operator liniar \c{s}i pozitiv, care reproduce func\c{t}ii constante, atunci pentru $f,g\in C[a,b]$ \c{s}i $x\in [a,b]$ fixat, au loc urm\u{a}toarele inegalit\u{a}\c{t}i:
		\begin{align}\label{inegCG2}
			\left|T(f,g;x)\right|&\leq \frac{1}{4}\cdot\widetilde{\omega}\left(f;2\cdot \sqrt{H(e_2;x)-H(e_1;x)^2}\right)\cdot\widetilde{\omega}\left(g;2\cdot \sqrt{H(e_2;x)-H(e_1;x)^2}\right)\\
			&\leq \frac{1}{4}\cdot \widetilde{\omega}\left(f;2\cdot \sqrt{H((e_1-x)^2;x)}\right)\cdot \widetilde{\omega}\left(g;2\cdot \sqrt{H((e_1-x)	 ^2;x)}\right)\nonumber.
		\end{align}
\end{theorem}

\begin{remark}
	Pentru a aplica inegalitatea de mai sus unora din cazurile operatorilor BEJ de tipul $I$ sau $II$, avem \^{\i}n primul r\^{a}nd nevoie de momentele de ordinul \^{\i}nt\^{a}i. Apoi trebuie s\u a calcul\u{a}m diferen\c{t}e de tipul $T(e_1,e_1;x):=H(e_2;x)-H(e_1;x)^2$. Dac\u a operatorul $H$ nu reproduce func\c tiile liniare, se ob\c tine o \^{i}mbun\u at\u a\c tire a inegalit\u a\c tii (\ref{inegCG1}).
\end{remark}

Operatorii Bernstein-Euler-Jacobi (BEJ) sunt \^{\i}mp\u{a}r\c{t}i\c{t}i \^{\i}n dou\u{a} clase: BEJ de tip $I$ \c{s}i BEJ de tip $II$. 
Scopul acestui articol este de a aplica inegalit\u{a}\c{t}ile de tip Chebyshev-Gr\"{u}ss de mai sus acestor tipuri de operatori. \^{I}n rezultatele urm\u{a}toare momentele de ordinul unu \c{s}i cele de ordinul doi vor fi de mare interes.  

Operatorii BEJ de tipul $I$, ca \c{s}i clas\u{a} de operatori pozitivi \c{s}i liniari, pot fi defini\c{t}i dup\u{a} cum urmeaz\u{a}. Pentru mai multe detalii, vezi \cite{Stanila:2012}. 

\begin{definition}
	Definim $R^{(r,a,b)}_{m,n}:C[0,1]\to C[0,1]$ ca fiind compozi\c{t}ia
		\begin{equation*}
			R^{(r,a,b)}_{m,n}=B_m \circ \mathcal{B}^{a,b}_r\circ B_n,
		\end{equation*}
	pentru $r>0$, $a,b \geq -1$, $n,m>1$. \^{I}n ecua\c{t}ia de mai sus, $\mathcal{B}^{a,b}_r$ este operatorul Euler-Jacobi-Beta definit \^{i}n \cite{Gonska-Rasa-Stanila:2014-2} \c{s}i 
	$B_n, B_m$ sunt operatorii Bernstein de ordin $n$ \c{s}i $m$. 
\end{definition}

\begin{remark}
	Operatorii  BEJ de tip $I$ reproduc constantele. Pentru anumite valori ale lui $a$, respectiv $b$, \c{s}i anume pentru $a=b=-1$, sunt reproduse \c{s}i func\c{t}iile liniare. 
\end{remark}

A doua clas\u{a} de operatori liniari \c{s}i pozitivi pe care \^{\i}i consider\u{a}m sunt BEJ de tipul $II$. Defini\c{t}ia este dat\u{a} \^{\i}n continuare.

\begin{definition}
	Pentru $r,s>0$, $a,b,c,d\geq -1$ \c{s}i $n>1$, definim $R^{s,c,d;r,a,b}_{n}:C[0,1]\to C[0,1]$ ca fiind
		\begin{equation*}
			R^{s,c,d;r,a,b}_n=\mathcal{B}^{c,d}_s\circ B_n\circ \mathcal{B}^{a,b}_r.
		\end{equation*}
	$\mathcal{B}^{a,b}_r$ \c{s}i $\mathcal{B}^{c,d}_s$ sunt operatori de tip Euler-Jacobi Beta \c{s}i $B_n$ este operatorul Bernstein de ordin $n$. 
\end{definition}

\begin{remark}
	Operatorii BEJ de tip $II$ reproduc constantele. Pentru anumite valori ale lui $a,b$ \c{s}i $c,d$, \c{s}i anume pentru $a=b=c=d=-1$, sunt reproduse \c{s}i func\c{t}iile liniare. 
\end{remark}

\section{Momentele de ordinul unu \c{s}i doi}

\begin{lemma}
Pentru clasa de operatori \textbf{BEJ de tipul I} momentele de ordinul 1 \c si 2 sunt date de:
\begin{equation}
		R^{(r,a,b)}_{m,n}((e_1-xe_0)^1;x)=\frac{a+1-x(a+b+2)}{r+a+b+2}. 
\end{equation}
\begin{equation}\label{bej1_mom}
\begin{array}{l}
R^{(r,a,b)}_{m,n}((e_1-xe_0)^2;x)=\dfrac{x^2[mn(a^2+b^2+5a+5b+2ab+6-r)+r^2(1-m-n)]}{mn(r+a+b+2)(r+a+b+3)}-\\ -\dfrac{x[mn(2a^2+2ab+8a+2b+6-r)+r^2(1-m-n)+mr(a-b)]}{mn(r+a+b+2)(r+a+b+3)}+\\ +\dfrac{[mn(a+1)(a+2)+m(r(a+1)+ab+a+b+1)]}{mn(r+a+b+2)(r+a+b+3)}
\end{array}
\end{equation}
\end{lemma}
Pentru detalii \^{\i}n ceea ce prive\c{s}te demonstra\c{t}ia, vezi \cite{Stanila:2012}.

\begin{lemma}
Pentru clasa de operatori \textbf{BEJ de tipul II} momentele de ordinul 1 \c si 2 sunt date de:
	\begin{align}
		R^{(s,c,d;r,a,b)}_{n}((e_1-xe_0)^1;x)=&\frac{-x\left[r(c+d+2)+(a+b+2)(s+c+d+2)\right]}{(r+a+b+2)(s+c+d+2)}+\\&+\frac{r(c+1)+(s+c+d+2)(a+1)}{(r+a+b+2)(s+c+d+2)}\nonumber. 
	\end{align}
	\begin{equation}\label{bej2_mom}
\begin{array}{l}
R^{(s,c,d;r,a,b)}_{n}((e_1-xe_0)^2;x)=\dfrac{(n-1)(a^2+b^2+2ab+5a+5b+6-r)(sx+c+1)(sx+c+2)}{n(r+a+b+2)(r+a+b+3)(s+c+d+2)(s+c+d+3)}+\\+\dfrac{(a^2+b^2+2ab+5a+5b+6-r-n(2ab+2a^2+8a+2b+6-r))(sx+c+1)}{n(r+a+b+2)(r+a+b+3)(s+c+d+2)}+\\+\dfrac{a^2+3a+2}{(r+a+b+2)(r+a+b+3)}+\dfrac{(sx+c+1)(sx+c+2)}{n(s+c+d+2)(s+c+d+3)}-\\-\dfrac{sx+c+1}{n(s+c+d+2)}-\dfrac{c^2+d^2+2cd+5c+5d+6-s}{(s+c+d+2)(s+c+d+3)}x^2+\\+\dfrac{2cd+2c^2+8c+2d+6-s}{(s+c+d+2)(s+c+d+3)}x-\dfrac{c^2+3c+2}{(s+c+d+2)(s+c+d+3)}+\\+\dfrac{2r(sx+c+1)}{n(r+a+b+2)(s+c+d+2)}-\dfrac{2r(sx+c+1)(sx+c+2)}{n(r+a+b+2)(s+c+d+2)(s+c+d+3)}+\\+\dfrac{2r(sx+c+1)(sx+c+2)}{(r+a+b+2)(s+c+d+2)(s+c+d+3)}+\dfrac{2(a+1-x(2r+a+b+2))(sx+c+1)}{(r+a+b+2)(s+c+d+2)}-\\-\dfrac{2(a+1)x}{r+a+b+2}+2x^2.
\end{array}
\end{equation}
\end{lemma}

Pentru detalii \^{\i}n ceea ce prive\c{s}te demonstra\c{t}ia, vezi \cite{Stanila:2012}.

Mul\c ti operatori uzuali pot fi reg\u asi\c ti dac\u a particulariz\u am valorile indicilor. 
Astfel \^in primele cinci tabele prezentate in Anex\u a reg\u asim momentele de ordinul doi, pentru cazurile particulare pe care am reu\c sit s\u a le localiz\u am \^in literatur\u a, determinate pe baza ecua\c tiilor (\ref{bej1_mom}) \c si (\ref{bej2_mom}), folosind  conven\c tia $B_{\infty} = \mathcal{B}_{\infty}^{a,b} = Id$.

\section{Inegalit\u a\c ti Chebyshev-Gr\"uss pentru operatorii BEJ de tip I \c si II}

\^{I}n continuare aplic\u am inegalitatea de tip Chebyshev-Gr\"{u}ss, dat\u a de (\ref{inegCG1}), operatorilor BEJ de tipul I \c si II \c si ob\c tinem urm\u atoarele teoreme.

\begin{theorem}
 Dac\u{a} $f,g\in C[a,b]$ \c{s}i $x\in [a,b]$ este fixat, atunci inegalitatea
 	\begin{equation*}
 		\left|T(f,g;x)\right|\leq \frac{1}{4}\widetilde{\omega}\left(f;2\cdot \sqrt{R^{(r,a,b)}_{m,n}((e_1-x)^2;x)}\right)\cdot \widetilde{\omega}\left(g;2\cdot \sqrt{R^{(r,a,b)}_{m,n}((e_1-x)^2;x)}\right)
 	\end{equation*}
 are loc, pentru $r> 0$, $a,b\geq -1$, $n,m>1$.  
\end{theorem}

\begin{theorem}
 Dac\u{a} $f,g\in C[a,b]$ \c{s}i $x\in [a,b]$ este fixat, atunci inegalitatea
 	\begin{equation*}
 		\left|T(f,g;x)\right|\leq \frac{1}{4}\widetilde{\omega}\left(f;2\cdot \sqrt{R^{(s,c,d;r,a,b)}_{n}((e_1-x)^2;x)}\right)\cdot \widetilde{\omega}\left(g;2\cdot \sqrt{R^{(s,c,d;r,a,b)}_{n}((e_1-x)^2;x)}\right)
 	\end{equation*}
 are loc, pentru $r,s>0$, $a,b,c,d\geq -1$, $n>1$.  
\end{theorem}

\begin{remark}
	Particulariz\^and valorile indicilor \^in ecua\c tiile (\ref{bej1_mom}) \c si (\ref{bej2_mom}), ob\c tinem inegalit\u a\c tile Chebyshev-Gr\" uss corespunz\u atoare operatorilor binecunoscu\c ti \^in literatur\u a. Momentele de ordinul doi \^in cazurile particulare sunt prezentate \^in Tabelele $1 - 5$ din Anex\u a. 
\end{remark}

Dac\u a aplic\u am inegalitatea (\ref{inegCG2}) operatorilor BEJ de tipul I \c si II care nu reproduc func\c tiile liniare, ob\c tinem urm\u atorul rezultat.

\begin{theorem}
 Pentru $f,g\in C[a,b]$ \c{s}i $x\in [a,b]$ fixat, urm\u atoarea inegalitate are loc:
 	\begin{equation*}
			\left|T(f,g;x)\right|\leq \frac{1}{4}\cdot\widetilde{\omega}\left(f;2\cdot \sqrt{T(e_1,e_1;x)}\right)\cdot\widetilde{\omega}\left(g;2\cdot \sqrt{T(e_1,e_1;x)}\right).
		\end{equation*}
 Diferen\c tele de mai sus sunt date de ecua\c tiile:
 	\begin{equation*}
 		T(e_1,e_1;x):= R^{(r,a,b)}_{m,n}((e_1-xe_0)^2;x)-[R^{(r,a,b)}_{m,n}((e_1-xe_0)^1;x)]^2,
 	\end{equation*}  
 respectiv
 	\begin{equation*}
 		T(e_1,e_1;x):= R^{(s,c,d;r,a,b)}_{n}((e_1-xe_0)^2;x)-[R^{(s,c,d;r,a,b)}_{n}((e_1-xe_0)^1;x)]^2.
 	\end{equation*}  
\end{theorem}

\begin{remark}
	Momentele de ordinul \^{i}nt\^{a}i \c si diferen\c tele de tip $T(e_1,e_1;x)$ pentru cazurile particulare se reg\u asesc \^{i}n Tabelele $6-7$ din Anex\u a. 
\end{remark}

\begin{landscape}
\section{Anex\u a}
\scalebox{0.9}{
\begin{tabular}{|C{2.5cm}|C{4cm}|C{13cm}|C{2cm}|}
\hline
Nota\c tie  &Denumire & Momentul de ordinul doi & Referin\c te\\
\hline
$U_n=R_{n,\infty}^{(n,-1,-1)}=R_{n}^{\infty,-,-;n,-1,-1}$& operatorul "original" Bernstein-Durrmeyer&  $U_n((e_1-xe_0)^2;x):=\dfrac{2X}{n+1}$ &\cite{Chen:1987} \cite{Goodman-Sharma:1991}\\\hline
$M_n=R_{n,\infty}^{(n,0,0)}=R_{n}^{\infty,-,-;n,0,0}$ & operatorul lui Durrmeyer &  $M_n((e_1-xe_0)^2;x):=\dfrac{2X(n-3)+2}{(n+2)(n+3)}$ &\cite{Durrmeyer:1986}\\
\hline
$D^{<\alpha>}=R_{n,\infty}^{(n, \alpha, \alpha)}=R_{n}^{\infty,-,-;n, \alpha, \alpha}$ &operatorul Bernstein-Durrmeyer cu ponderi simetrice&  $D^{<\alpha>}((e_1-xe_0)^2;x):=\dfrac{X(2n-4a^2-10a-6)+(a+1)(a+2)}{(n+2a+2)(n+2a+3)}$& \cite{Lupas:1995}\\
\hline
$M_n^{ab}=R_{n,\infty}^{(n,a,b)}=R_{n}^{\infty,-,-;n,a,b}$ & operatorul lui Durrmeyer cu ponderi Jacobi& $M_n^{ab}((e_1-xe_0)^2;x):=\dfrac{a^2+b^2+2ab+5a+5b+6-2n}{(n+a+b+2)(n+a+b+3)}\cdot x^2$&\\ &&$-\dfrac{2a^2+2ab+8a+2b+6-2n}{(n+a+b+2)(n+a+b+3)}\cdot x +\dfrac{(a+1)(a+2)}{(n+a+b+2)(n+a+b+3)}$& \cite{Paltanea:1983}\\
\hline
$U_n^{\varrho}=R_{n,\infty}^{(n\varrho,-1,-1)}=R_{n}^{\infty,-,-;n\varrho,-1,-1}$ &$\varrho$-operatorul lui P\u alt\u anea& $U_n^{\varrho}((e_1-xe_0)^2;x):=\dfrac{X(1+\varrho)}{n\varrho+1}$& \cite{Paltanea:2007}\\
\hline
$P_n=R_{n,\infty}^{(nc,a,b)}=R_{n}^{\infty,-,-;nc,a,b}$ &operatorul lui Mache-Zhou& $P_n((e_1-xe_0)^2;x):=\dfrac{a^2+b^2+2ab+5a+5b+6-nc-nc^2}{(nc+a+b+2)(nc+a+b+3)}\cdot x^2$&\\ &&$-\dfrac{(2a^2+2ab+8a+2b+6-nc-nc^2)\cdot x+(a+1)(a+2)}{(nc+a+b+2)(nc+a+b+3)}$&\cite{Mache-Zhou:1996}\\
\hline
$ L_n^{\nabla}=R_{\infty,n}^{(n,-1,-1)}=R_{n}^{n,-1,-1;\infty,-,-}$ &operator de tip Stancu &$ L_n^{\nabla}((e_1-xe_0)^2;x):=\dfrac{2X}{2+1}$ &\cite{Lupas-Lupas:1987}\\
\hline
$V_n^{0,0}=R_{\infty,n}^{(n,0,0)}=R_{n}^{n,0,0;\infty,-,-}$& operatorul lui Lupa\c s & $V_n^{0,0}((e_1-xe_0)^2;x):=\dfrac{X(2n^2-6n)+3n+1}{n(n+2)(n+3)}$& \cite{Lupas:1995}\\
\hline

\end{tabular}}
\captionof*{table}{\textbf{Tabelul 1}: Momentele de ordinul doi - partea I}
\scalebox{0.9}{
\begin{tabular}{|C{2.9cm}|C{4cm}|C{13cm}|C{2cm}|}
\hline
Nota\c tie  &Denumire & Momentul de ordinul doi & Referin\c te\\
\hline
$V_n^{\alpha,\beta}=R_{\infty,n}^{(n, \alpha, \beta)}=R_{n}^{n, \alpha, \beta;\infty,-,-}$& operatorul lui Lupa\c s cu ponderi Jacobi& $V_n^{\alpha,\beta}((e_1-xe_0)^2;x):=\dfrac{\alpha^2+\beta^2+2\alpha\beta +5\alpha+5\beta+6-2n}{(n+\alpha+\beta+2)(n+\alpha+\beta+3)}\cdot x^2$&\\ &&$-\dfrac{2\alpha^2+2\alpha\beta+9\alpha+\beta+6-2n}{(n+\alpha+\beta+2)(n+\alpha+\beta+3)}\cdot x +\dfrac{n\alpha^2+4n\alpha+\alpha\beta+3n+\alpha+\beta+1}{n(n+\alpha+\beta+2)(n+\alpha+\beta+3)}$& I. Ra\c sa handwritten notes, 19 August 2008.\\
\hline
$S_n^{\alpha}=R_{\infty,n}^{(\frac{1}{\alpha}, -1, -1)}=R_{n}^{\frac{1}{\alpha}, -1, -1;\infty,-,-}$& operatorul lui Stancu& $S_n^{\alpha}((e_1-xe_0)^2;x):=\dfrac{X(n^2+1)}{n(n+1)}$&\cite{Stancu:1968}\\
\hline
$Q_n^{\varrho,c,d}=R_{\infty,n \varrho}^{(n,c,d)}=R_{n}^{n \varrho, c, d;\infty,-,-}$ & operator de tip Stancu cu parametri $c$ \c si $d$& $Q_n^{\varrho,c,d}((e_1-xe_0)^2;x):=\dfrac{c^2+d^2+2cd+5c+5d+6-n\varrho-n\varrho^2}{(n\varrho+c+d+2)(n\varrho+c+d+3)}\cdot x^2$&\\ &&$-\dfrac{2c^2+2cd+8c+2d+6-n\varrho-n\varrho^2}{(n\varrho+c+d+2)(n\varrho+c+d+3)}\cdot x +\dfrac{nc^2+3nc+nc\varrho+n\varrho+2n+c+d+cd+1}{n(n\varrho+c+d+2)(n\varrho+c+d+3)}$&H. Gonska handwritten notes, 18 March 2009.\\
\hline
$\mathbb{\overline{B}}_n$, $\mathcal{B}_n^{-1,-1}=R_{\infty,\infty}^{(n,-1,-1)}=R_{\infty}^{\infty,-,-;n,-1,-1}=R_{\infty}^{n,-1,-1;\infty,-,-}$& operatorul Beta de a doua spe\c ta al lui Lupa\c s  &$\mathbb{\overline{B}}_n((e_1-xe_0)^2;x):=\dfrac{X}{n+1}$ &\cite{Lupas:1972}\\
\hline
$\mathbb{\widetilde{B}}_{\lambda}$, $T_{\lambda}=R_{\infty,\infty}^{(\frac{1}{\lambda},-1,-1)}=R_{\infty}^{\infty,-,-;\frac{1}{\lambda},-1,-1}=R_{\infty}^{\frac{1}{\lambda},-1,-1;\infty,-,-}$&operatorul Beta al lui M\"uhlbach &$\mathbb{\widetilde{B}}_{\lambda}((e_1-xe_0)^2;x):=\dfrac{\lambda X}{1+\lambda}$ &\cite{Muhlbach:1972}\\
\hline
$\mathbb{B}_n$, $\mathcal{B}_n^{0,0}=R_{\infty,\infty}^{(n,0,0)}=R_{\infty}^{\infty,-,-;n,0,0}=R_{\infty}^{n,0,0;\infty,-,-}$& operatorul Beta de prima spe\c ta al lui Lupa\c s &$\mathbb{B}_n((e_1-xe_0)^2;x):=\dfrac{X(n-6)+2}{(n+2)(n+3)}$ &\cite{Lupas:1972}\\
\hline
\end{tabular}}
\captionof*{table}{\textbf{Tabelul 2}: Momentele de ordinul doi - partea a II-a}
\scalebox{0.9}{
\begin{tabular}{|C{4cm}|C{3cm}|C{12cm}|C{2cm}|}
\hline
Nota\c tie  &Denumire & Momentul de ordinul doi & Referin\c te\\
\hline
$\mathcal{B}_n^{-1,\beta}=R_{\infty,\infty}^{(n,-1,\beta)}=R_{\infty}^{\infty,-,-;n,-1, \beta}=R_{\infty}^{n,-1, \beta;\infty,-,-}$& &$\mathcal{B}_n^{-1,\beta}((e_1-xe_0)^2;x):=\dfrac{nX+(\beta+1)(\beta+2)x^2}{(n+\beta+1)(n+\beta+2)}$ &\cite{Gonska-Rasa-Stanila:2014-2}\\
\hline
$\mathcal{B}_n^{\alpha,-1}=R_{\infty,\infty}^{(n,\alpha,-1)}=R_{\infty}^{\infty,-,-;n, \alpha,-1}=R_{\infty}^{n, \alpha,-1;\infty,-,-}$& &$\mathcal{B}_n^{\alpha,-1}((e_1-xe_0)^2;x):=\dfrac{nX+(\alpha+1)(\alpha+2)(x-1)^2}{(n+\alpha+1)(n+\alpha+2)}$ &\cite{Gonska-Rasa-Stanila:2014-2}\\
\hline
$\mathcal{B}_n^{\alpha,\beta}=R_{\infty,\infty}^{(n,\alpha, \beta)}=R_{\infty}^{\infty,-,-;n, \alpha, \beta}=R_{\infty}^{n, \alpha, \beta;\infty,-,-}$& operatorul Beta cu ponderi Jacobi &$\mathcal{B}_n^{\alpha,\beta}((e_1-xe_0)^2;x):=\dfrac{\alpha^2+\beta^2+2\alpha\beta+5\alpha+5\beta+6-n}{(n+\alpha+\beta+2)(n+\alpha+\beta+3)}\cdot x^2-\dfrac{2\alpha^2+2\alpha\beta+8\alpha+2\beta+6-n}{(n+\alpha+\beta+2)(n+\alpha+\beta+3)}\cdot x+\dfrac{(\alpha+1)(\alpha+2)}{(n+\alpha+\beta+2)(n+\alpha+\beta+3)}$ &\cite{Gonska-Rasa-Stanila:2014-2}\\
\hline
$B_n=R_{n,\infty}^{(\infty,-, -)}=R_{\infty,n}^{(\infty,-, -)}=R_{n}^{\infty,-,-;\infty,-,-}$& operatorul lui Bernstein & $B_n((e_1-xe_0)^2;x):=\dfrac{X}{n}$&\cite{Bernstein:1912}\\
\hline
\end{tabular}}
\captionof*{table}{\textbf{Tabelul 3}: Momentele de ordinul doi - partea a III-a}

\scalebox{0.9}{
\begin{tabular}{|C{4cm}|C{3cm}|C{12cm}|C{2cm}|}
\hline
Nota\c tie  &Denumire & Momentul de ordinul doi & Referin\c te\\
\hline
$D_n=R_{n,n +1}^{(\infty,-,-)}$& &$D_n((e_1-xe_0)^2;x):=\dfrac{2X}{n+1}$ &\cite{Gonska-Pitul-Rasa:2006-1}\\
\hline
$R_{m,n}^{\infty}=R_{m,n}^{(\infty,-,-)}$&&$R_{m,n}^{\infty}((e_1-xe_0)^2;x):=\dfrac{X(n+m-1)}{nm}$&\cite{Stanila:2012}\\
\hline
$R_{m,n}^{\varrho}=R_{m,n}^{(n\varrho,-1,-1)}$&&$R_{m,n}^{\varrho}((e_1-xe_0)^2;x):=\dfrac{X(n\varrho+m\varrho-\varrho+m)}{m(n\varrho+1)}$&\cite{Stanila:2012}\\
\hline
\end{tabular}}
\captionof*{table}{\textbf{Tabelul 4}: Momentele de ordinul doi - pentru operatori BEJ de primul tip}
\newpage

\scalebox{0.9}{
\begin{tabular}{|C{4cm}|C{3cm}|C{12cm}|C{2cm}|}
\hline
Nota\c tie  &Denumire & Momentul de ordinul doi & Referin\c te\\
\hline
$\mathbb{B}_{\infty}^{(\alpha,\lambda)}=R_{\infty}^{\frac{1}{\alpha}, -1, -1;\frac{1}{\lambda}, -1, -1}$& operatori de tip Beta generaliza\c ti&$\mathbb{B}_{\infty}^{\alpha,\lambda}((e_1-xe_0)^2;x):=\dfrac{X(\alpha+\lambda+\alpha\lambda)}{(1+\alpha)(1+\lambda)}$ &\cite{Pitul:2007}\\
\hline
$F_{n}^{\alpha}=R_{n}^{\frac{1}{\alpha}, -1, -1;n -1, -1}$&operatorul lui Finta &$F_{n}^{\alpha}((e_1-xe_0)^2;x):=\dfrac{X(n\alpha+\alpha+2)}{(\alpha+1)(n+1)}$&\cite{Finta:2002}\\
\hline
$\mathbb{B}_{n}^{(\alpha,\lambda)}=R_{n}^{\frac{1}{\alpha}, -1, -1;\frac{1}{\lambda}, -1, -1}$&operatorul de tip Beta al lui Pi\c tul&$\mathbb{B}_{n}^{(\alpha,\lambda)}((e_1-xe_0)^2;x):=\dfrac{X(n\alpha+n\lambda+n\alpha\lambda+1)}{n(1+\alpha)(1+\lambda)}$&\cite{Pitul:2007}\\
\hline
\end{tabular}}
\captionof*{table}{\textbf{Tabelul 5}: Momentele de ordinul doi - pentru operatori BEJ de tipul al doilea}

\scalebox{0.9}{
\begin{tabular}{|C{1.5cm}|C{4cm}|C{17cm}|}
\hline
Nota\c tie  &Momentul de ordinul unu & Diferen\c te de tipul $T(e_1,e_1;x)$\\
\hline
$M_n$ & $M_n((e_1-xe_0)^1;x):=\dfrac{1-2\cdot x}{n+2}$ & $T(e_1,e_1;x):=\dfrac{(2\cdot n\cdot X+1)\cdot (n+1)}{(n+2)^2\cdot (n+3)}$\\
\hline
$D^{<\alpha>}$ & $D^{<\alpha>}((e_1-xe_0)^1;x):=\dfrac{-2\cdot x\cdot (\alpha+1)+\alpha+1}{n+2\cdot \alpha+2}$ &  $T(e_1,e_1;x):=\dfrac{X\cdot (2\cdot n^2+2\cdot n+2\cdot \alpha\cdot n)+1+n+2\alpha +\alpha n+\alpha^2}{(n+2\alpha+2)^2 (n+2\alpha+3)}$\\
\hline
$M_n^{ab}$ & $M^{ab}_n((e_1-xe_0)^1;x):=\dfrac{a+1-x\cdot (a+b+2)}{n+a+b+2}$& $T(e_1,e_1;x):=\dfrac{2nX(n+b+1)+n\cdot x^2(b-a)+1+n+a+b+an+ab}{(n+a+b+2)^2(n+a+b+3)}$\\
\hline
$P_n$ &$P_n((e_1-xe_0)^1;x):=\dfrac{a+1-x(a+b+2)}{n\cdot c+a+b+2}$& $T(e_1,e_1;x):=\dfrac{X\cdot n\cdot c^2(n\cdot c+a+b+n+2)+x\cdot n\cdot c(a+b)}{(nc+a+b+2)^2(nc+a+b+3)}$\\ &&$-\dfrac{10a^2+5anc+15a+3b+2a^2nc+3nc+2a^2b+2a^3+5ab+7}{(nc+a+b+2)^2(nc+a+b+3)}$\\
\hline
$V_n^{0,0}$& $V^{0,0}_n((e_1-xe_0)^1;x):=\dfrac{1-2x}{n+2}$ & $T(e_1,e_1;x):=\dfrac{2(n+1)(n^2\cdot X+n+1)}{n(n+2)^2 (n+3)}$\\
\hline
$V_n^{\alpha,\beta}$& $V^{\alpha,\beta}_n((e_1-xe_0)^1;x):=\dfrac{\alpha+1-x\cdot (\alpha+\beta+2)}{n+\alpha+\beta+2}$& $T(e_1,e_1;x):=\dfrac{2n^3+n^2\beta+n^2\alpha+2n^2}{(n+\alpha+\beta+2)^2(n+\alpha+\beta+3)n}\cdot x^2$\\ &&$+\dfrac{n\alpha^2-2n^2+n^2\alpha+2n\alpha -n\beta^2 -2\beta n-2n^3 -3n^2\beta}{(n+\alpha+\beta+2)^2(n+\alpha+\beta+3)n}\cdot x +\dfrac{-2-4n-3\alpha-3\beta-n\alpha^2-4\alpha\beta -3\alpha\beta n-\alpha^2-\beta^2 -5n\alpha -\beta^2\alpha-3\beta n-2n^2-\alpha^2\beta-2n^2\alpha}{n(n+\alpha+\beta+2)^2(n+\alpha+\beta+3)}$\\
\hline
\end{tabular}}
\captionof*{table}{\textbf{Tabelul 6}: Diferen\c te de tipul $T(e_1,e_1;x)$- Partea I}
\scalebox{0.9}{
\begin{tabular}{|C{2cm}|C{6cm}|C{14cm}|}
\hline
Nota\c tie  &Momentul de ordinul unu & Diferen\c te de tipul $T(e_1,e_1;x)$ \\
\hline
$Q_n^{\rho,c,d}$& $Q^{\rho,c,d}_n((e_1-xe_0)^1;x):=\dfrac{c+1-x\cdot (c+d+2)}{n\cdot \rho+c+d+2}$& $T(e_1,e_1;x):=\dfrac{2n^2\rho^2+n^3\rho^2+n^2\rho^2 d+n^2c\rho^2+n^3\rho^3}{(n\rho+c+d+2)^2(n\rho+c+d+3)n}\cdot x^2$\\ &&$-\dfrac{-n^2 c\rho+n^2 d \rho+2n^2\rho^2+n^3\rho^2 +n^2\rho^2 d +n^2 c \rho^2+n^3\rho^3}{(n\rho+c+d+2)^2 (n\rho+c+d+3) n}\cdot x +\dfrac{-2-n-4nc\rho-2cdn\rho-n^2\rho-nd -3\rho n-3c-nc -4cd -c^2-d^2-c^2n\rho}{n(n\rho+c+d+2)^2(n\rho+c+d+3)}$\\&&$+\dfrac{-n^2 c\rho-n^2 c \rho^2-ncd-3d-d^2c-n^2\rho^2-c^2d-2dn\rho}{n(n\rho+c+d+2)^2(n\rho+c+d+3)}$\\
\hline
$\mathbb{B}_n$, $\mathcal{B}_n^{0,0}$& $\mathbb{B}_n((e_1-xe_0)^1;x):=\dfrac{1-2x}{n+2}$ &$T(e_1,e_1;x):=\dfrac{X(n^2+12n+24)+n+1}{(n+2)^2(n+3)}$ \\
\hline
$\mathcal{B}_n^{-1,\beta}$& $\mathcal{B}_n^{-1,\beta}((e_1-xe_0)^1;x):=\dfrac{-x(\beta+1)}{n+\beta+1}$ & $T(e_1,e_1;x):=\dfrac{n^2\cdot X+n\beta+n}{(n+\beta+1)^2(n+\beta+2)}$ \\
\hline
$\mathcal{B}_n^{\alpha,-1}$& $\mathcal{B}_n^{\alpha,-1}((e_1-xe_0)^1;x):=\dfrac{(\alpha+1)(1-x)}{n+\alpha+1}$ & $T(e_1,e_1;x):=\dfrac{n^2\cdot X-n\cdot x (1+\alpha)+n+n\cdot \alpha}{(n+\alpha+1)^2 (n+\alpha+2)}$\\
\hline
$\mathcal{B}_n^{\alpha,\beta}$& $\mathcal{B}_n^{\alpha,\beta}((e_1-xe_0)^1;x):=\dfrac{\alpha+1-x(\alpha+\beta+2)}{n+\alpha+\beta+2}$ & $T(e_1,e_1;x):=\dfrac{X\cdot n^2+nx(\beta-\alpha)+\alpha+\beta+\alpha\beta+\alpha n+n+1}{(n+\alpha+\beta+2)^2 (n+\alpha+\beta+3)}$ \\
\hline
\end{tabular}}
\captionof*{table}{\textbf{Tabelul 7}: Diferen\c te de tipul $T(e_1,e_1;x)$- Partea a II-a}
\end{landscape}

\begingroup
\renewcommand{\section}[2]{{\Large\bf Bibliografie}}%

\endgroup
\end{document}